\input pictex.tex
\input amstex.tex
\documentstyle{amsppt}
\magnification=1200
 \pagewidth{16.3truecm}
 \pageheight{24truecm}
 \nologo

\def\div{\operatorname{div}}
\def\mod{\operatorname{mod}}
\def\Spec{\operatorname{Spec}}
\def\supp{\operatorname{supp}}
\def\Var{\operatorname{Var}}
\def\ord{\operatorname{ord}}

\def\Var{\operatorname{Var}}
\def\dim{\operatorname{dim}}

\refstyle{A} \widestnumber\key{AKMW}
 \topmatter
\title
On motivic principal value integrals
\endtitle
\author
Willem Veys  
\endauthor
\address K.U.Leuven, Departement Wiskunde, Celestijnenlaan 200B,
         B--3001 Leuven, Belgium  \endaddress
\email wim.veys\@wis.kuleuven.ac.be  \newline
 http://www.wis.kuleuven.ac.be/algebra/veys.htm
\endemail
 \abstract Inspired by $p$-adic (and real) principal value
 integrals, we introduce motivic principal value integrals
 associated to multi-valued rational differential forms on smooth algebraic varieties.
 We investigate the natural question whether (for complete varieties) this notion
 is a birational invariant. The answer seems to be related to the
 dichotomy of the Minimal Model Program.
\endabstract
 \endtopmatter

\document
\noindent {\bf Introduction}
 \bigskip
 \noindent {\bf 0.1.} Real
and $p$-adic principal value integrals were first introduced by
Langlands in the study of orbital integrals [Lan1] [Lan2] [LS1]
[LS2]. They are associated to multi-valued differential forms on
real and $p$-adic manifolds, respectively.

Let for instance $X$ be a complete smooth connected algebraic
variety of dimension $n$ over $\Bbb Q_p$ (the field of $p$-adic
numbers). Denoting by $\Omega^n_X$ the vector space of {\it
rational} differential $n$-forms on $X$, take $\omega \in
(\Omega^n_X)^{\otimes d}$ defined over $\Bbb Q_p$; we then write
formally $\omega^{1/d}$ and consider it as a multi-valued rational
differential form on $X$.

We suppose that $\div \omega$ is a normal crossings divisor (over
$\Bbb Q_p$) on $X$; say $E_i, i \in S$,  \linebreak   are its
irreducible components. Let $\div \omega^{1/d} := \frac 1 d \div
\omega = \sum_{i \in S} (\alpha_i - 1) E_i$, where then the
$\alpha_i \in \frac 1 d \Bbb Z$. If $\omega^{1/d}$ has {\it no
logarithmic poles}, i.e. if all $\alpha_i \ne 0$, the {\it
principal value integral} $PV \int_{X(\Bbb Q_p)} |
\omega^{1/d}|_p$ of $\omega^{1/d}$ on $X(\Bbb Q_p)$ is defined as
follows.  Cover $X(\Bbb Q_p)$ by (disjoint) small enough open
balls $B$ on which there exist local coordinates $x_1, \dots ,
x_n$ such that all $E_i$ are coordinate hyperplanes. Consider for
each $B$ the {\sl converging} integral $\int_B |x_1 x_2 \cdots x_n
|^s_p | \omega^{1/d}|_p$ for $s \in \Bbb C$ with ${\Cal R}(s) >\!>
0$, take its meromorphic continuation to $\Bbb C$ and evaluate
this in $s = 0$; then add all these contributions. One can check
that the result is independent of all choices.

In the real setting we proceed similarly but then we also need a
partition of unity, and we have to assume that $\omega^{1/d}$ has
no integral poles, i.e. the $\alpha_i \notin \Bbb Z_{\leq 0}$.
Here the independency result is somewhat more complicated; it was
verified in detail in [Ja1].

\bigskip
\noindent {\bf 0.2.} These principal value integrals appear as
coefficients of asymptotic expansions of oscillating integrals and
fibre integrals, and as residues of poles of distributions
$|f|^\lambda$ or Igusa zeta functions. See [Ja1, \S 1] for an
overview and [AVG][De2][Ig1][Ig2][Ja2][Lae] for more details.

\bigskip
\noindent {\bf 0.3.} Last years (usual) $p$-adic integration and
$p$-adic Igusa zeta functions were `upgraded' to motivic
integration and motivic zeta functions in various important papers
of Denef and Loeser (after an idea of Kontsevich [Ko]). We mention
the first papers [DL1] [DL2] and surveys [DL3][Lo][Ve4].

In this note we introduce similarly motivic principal value
integrals. It is not totally clear what the most natural approach
is; however the following should be satisfied. Returning to the
$p$-adic setting of (0.1), we denote $E^\circ_I := (\cap_{i \in
I}E_i) \setminus (\cup_{\ell \not\in I} E_\ell)$ for $I \subset
S$. So $X = \coprod_{I \subset S} E^\circ_I$. Then, if suitable
conditions about good reduction $\mod p$ are satisfied, a similar
proof as for Denef's formula for the $p$-adic Igusa zeta function
[De1] yields that $PV \int_{X(\Bbb Q_p)} |\omega^{1/d}|_p$ is
given (up to a power of $p$) by
$$\sum_{I \subset S} \sharp (E^\circ_I)_{\Bbb F_p} \prod_{i \in I}
\frac{p-1}{p^{\alpha_i}-1},$$ where $\sharp(\cdot)_{\Bbb F_p}$
denotes the number of $\Bbb F_p$-rational points of the reduction
$\mod p$. Since motivic objects should specialize to the analogous
$p$-adic objects (for almost all $p$), any decent definition of a
motivic principal value integral $PV\int_X \omega^{1/d}$
associated to analogous $X$ and $\omega^{1/d}$ (say over $\Bbb C$)
should boil down to the formula
$$\sum_{I \subset S} [E^\circ_I] \prod_{i \in I}
\frac{L-1}{L^{\alpha_i}-1}$$ (up to a power of $L$). Here $[ \cdot
]$ denotes the class of a variety in the Grothendieck ring of
algebraic varieties, and $L := [ \Bbb A^1]$, see (1.5). Note also
that this is precisely the `user-friendly formula' (in the
terminology of [Cr]) for the {\sl converging} motivic integral
associated to the $\Bbb Q$-divisor $\div \omega^{1/d} = \sum_{i
\in S} (\alpha_i - 1)E_i$ {\sl if all} $\alpha_i > 0$. We will use
(evaluations of) motivic zeta functions as in [Ve2] or [Ve3] to
introduce this desired motivic principal value integral.

\bigskip
\noindent {\bf 0.4.} {\sl Remark.} (1) Remembering the origin of
principal value integrals, we mention that Hales introduced
motivic orbital integrals, specializing to the usual $p$-adic
orbital integrals [Ha1][Ha2].

(2) As in the $p$-adic case, the study of motivic principal value
integrals, especially their vanishing, is related to determining
the poles of motivic zeta functions, and hence of the derived
Hodge and topological zeta functions. A nice result about the
vanishing of real principal value integrals, and a conjecture in
the $p$-adic case, is in [DJ].

\bigskip
\noindent {\bf 0.5.} Since a multi-valued differential form
$\omega^{1/d}$ is in fact a birational notion, it is a natural
question whether the motivic principal value integral is a
birational invariant. In other words, if $X_1$ and $X_2$ are
different complete smooth models of the birational equivalence
class associated to $\omega^{1/d}$ such that $\div \omega^{1/d}$
is a normal crossings divisor and $\omega^{1/d}$ has no
logarithmic poles on both $X_1$ and $X_2$, is then $PV \int_{X_1}
\omega^{1/d} = PV \int_{X_2} \omega^{1/d}$ ?

This appears to be related to the dichotomy of the Minimal Model
Program. We show by explicit counterexamples that the answer is in
general negative when the Kodaira dimension is $-\infty$. On the
other hand, when the Kodaira dimension is nonnegative, we prove
birational invariance in dimension $2$. In higher dimensions, we
explain how the motivic principal value integrals yield a
`partial' birational invariant, assuming the Minimal Model
Program. Here some subtle problems appear, which we think are
interesting to investigate.

\bigskip
\noindent {\bf 0.6.} We also introduce motivic principal value
integrals on a smooth variety $X$ if $\div \omega^{1/d}$ is not
necessarily a normal crossings divisor, facing similar problems.
(For real principal value integrals this was considered by Jacobs
[Ja1, \S 7].)

\bigskip
\noindent {\bf 0.7.} We will work over a field $k$ of
characteristic zero. When using minimal models we moreover assume
$k$ to be algebraically closed. In \S1 we briefly recall the
necessary birational geometry, and the for our purposes relevant
motivic zeta function. In \S 2 we introduce motivic principal
value integrals on smooth varieties. We first proceed on the level
of Hodge polynomials to show that our approach with evaluations of
motivic zeta functions can really be considered as the analogue of
`classical' real or $p$-adic principal value integrals. Then in \S
3 we consider the birational invariance question.

\bigskip
\bigskip
\noindent {\bf 1. Birational geometry}
\bigskip
\noindent As general references for \S \S 1.1--1.4 we mention [KM]
and [Ma].
\bigskip
\noindent {\bf 1.1.} An algebraic variety is an integral separated
scheme of finite type over $\Spec k$, where $k$ is a field of
characteristic zero. A modification is a proper birational
morphism. A log resolution of an algebraic variety is a
modification $h : Y \rightarrow X$ from a smooth $Y$ such that the
exceptional locus of $h$ is a (simple) normal crossings divisor.

Let $B$ be a $\Bbb Q$-divisor on $X$. Then a log resolution of $B$
is a modification $h : Y \rightarrow X$ from a smooth $Y$ such
that the exceptional locus of $h$ is a divisor, and its union with
$h^{-1} (\supp B)$ is a (simple) normal crossings divisor.

\bigskip
\noindent {\bf 1.2.} Moreover, let $X$ be normal and denote $n :=
\dim X$. A (Weil) $\Bbb Q$-divisor $D$ on $X$ is called $\Bbb
Q$-Cartier if some integer multiple of $D$ is Cartier. And $X$ is
called $\Bbb Q$-factorial if every Weil divisor on $X$ is $\Bbb
Q$-Cartier.

The variety $X$ has a well-defined linear equivalence class $K_X$
of canonical (Weil) divisors. Its representatives are the divisors
$\div \eta$ of rational differential $n$-forms $\eta$ on $X$.
Denoting by $\Omega^n_X$ the vector space of those rational
differential $n$-forms, we can consider more generally elements
$\omega \in (\Omega^n_X)^{\otimes d}$ for any $d \in \Bbb Z_{>
0}$, and their associated divisor $\div \omega$. Then we write
formally $\omega^{1/d}$, considered as a multi-valued rational
differential form on $X$, and put $\div \omega^{1/d} := \frac 1d
\div \omega$.  Since $\div \omega$ represents $dK_X$, we can say
that the $\Bbb Q$-divisor $\div \omega^{1/d}$ represents $K_X$.

One says that $X$ is Gorenstein if $K_X$ is Cartier, and $\Bbb
Q$-Gorenstein if $K_X$ is $\Bbb Q$-Cartier.

\bigskip
\noindent {\bf 1.3.} For a $\Bbb Q$-Gorenstein $X$, let $h : Y
\rightarrow X$ be a log resolution of $X$, and denote by $E_i, i
\in S$, the irreducible components of the exceptional locus of
$h$. One says that $X$ is {\it terminal} and {\it canonical} if in
the expression
$$K_{Y|X} := K_Y - h^\ast K_X = \sum_{i \in S} a_i E_i$$
all $a_i, i \in S$, are greater than $0$ and at least $0$,
respectively. (These notions are independent of the chosen
resolution.) Such varieties can be considered `mildly' singular;
note that a smooth variety is terminal.

\bigskip
\noindent {\bf 1.4.} Let $k$ be algebraically closed.

 (i) A {\it minimal model} in a given
birational equivalence class of nonnegative Kodaira dimension is a
complete variety $X_m$ in this class which is $\Bbb Q$-factorial
and terminal and such that $K_{X_m}$ is nef. This last condition
means that  the intersection number $K_{X_m} \cdot C \geq 0$ for
all irreducible curves $C$ on $X_m$.

The existence of these objects is the heart of Mori's Minimal
Model Program. This is now accomplished in dimension $\leq 3$ (and
there is a lot of progress in dimension 4). In dimension 2 it is
well known that there is a unique minimal model, which is moreover
smooth, in each birational equivalence class. Also, each smooth
complete surface in the class maps to the unique minimal model
through a sequence of blowing-ups. In higher dimensions, two
different minimal models are isomorphic in codimension one. Here
each smooth complete variety in the class maps to a minimal model
through a {\it rational map} (which is a composition of divisorial
contractions and flips).

(ii) In a given birational equivalence class of general type (i.e.
of maximal Kodaira dimension), a {\it canonical model} is a
complete variety $X_c$ in this class which is canonical and such
that $K_{X_c}$ is ample. This object is unique and there is a
morphism from every minimal model in the class to it.

\bigskip
\noindent {\bf 1.5.}   Here, by abuse of terminology, we allow a
variety to be reducible.

(i) The Grothendieck ring $K_0 (Var_k)$ of algebraic varieties
over $k$ is the free abelian group generated by the symbols $[V]$,
where $[V]$ is a variety, subject to the relations $[V] =
[V^\prime]$ if $V$ is isomorphic to $V^\prime$, and $[V] = [V
\setminus W] + [W]$ if $W$ is closed in $V$. Its ring structure is
given by $[V] \cdot [W] := [V \times W]$. (See [Bi] for
alternative descriptions of $K_0(Var_k)$, and see [Po] for the
recent proof that it is not a domain.) Usually, one abbreviates $L
:= [\Bbb A^1]$.

For the sequel we need to extend $K_0(\Var_k)$ with fractional
powers of $L$ and to localize. Fix $d \in \Bbb Z_{> 0}$; we
consider
$$K_0 (\Var_k)[L^{-1/d}] := {K_0(\Var_k)[T] \over (LT^d - 1)}$$
(where $L^{-1/d} := \bar T$). We then localize this ring with
respect to the elements $L^{i/d} - 1$,  \linebreak
  $i \in \Bbb Z
\setminus \{ 0 \}$. What we really need is the subring of this
localization generated by $K_0 (\Var_k)$, $L^{-1}$ and the
elements $(L-1)/(L^{i/d} - 1), i \in \Bbb Z \setminus \{ 0 \}$; we
denote this subring by ${\Cal R}_d$.

(ii) For a variety $V$, we denote by $h^{p,q} (H^i_c (V, \Bbb C))$
the rank of the $(p,q)$-Hodge component in the mixed Hodge
structure of the $i$th cohomology group with compact support of
$V$. The {\it Hodge polynomial} of $V$ is
$$H(V) = H(V;u,v) := \sum_{p,q} \big(\sum_{i \geq 0} (-1)^i h^{p,q}
(H^i_c (V,\Bbb C))\big) u^pv^q \in \Bbb Z [u,v].$$ Precisely by
the defining relations of $K_0 (\Var_k)$, there is a well-defined
ring homomorphism $H : K_0 (\Var_k) \rightarrow \Bbb Z[u,v]$,
determined by $[V] \mapsto H(V)$. It induces a ring homomorphism
$H$ from ${\Cal R}$ to the `rational functions in $u,v$ with
fractional powers'.

\bigskip
\noindent {\bf 1.6.} In [DL1] Denef and Loeser associated a
motivic zeta function to a regular function on a smooth variety.
In [Ve2] and [Ve3, \S 2] we considered several generalizations; we
mention here a special case of [Ve3, (2.2)].

(i) Let $X$ be a canonical variety and $D$ any $\Bbb Q$-divisor on
$X$. Take a log resolution $h : Y \rightarrow X$ of $D$ and denote
by $E_i, i \in S$, the irreducible components of the union of
$h^{-1} (\supp D)$ and the exceptional locus of $h$. For each $i
\in S$ let $\nu_i - 1$ and $N_i$ denote the multiplicity of $E_i$
in $K_{Y|X}$ and $h^\ast D$, respectively. Note that all $\nu_i
\geq 1$ since $X$ is canonical. We also put $E^\circ_I := (\cap_{i
\in I} E_i) \setminus (\cup_{\ell \not\in I} E_\ell)$ for $I
\subset S$. We associated to $D$ on $X$ the zeta function
 $${\Cal Z}_X(D;s) := L^{-n} \sum_{I \subset S} [E^\circ_I]
\prod_{i \in I} \frac{L-1}{L^{\nu_i + sN_i}-1}.$$
 Here $L^{-s}$ is just the traditional notation  for a variable $T$. So
${\Cal Z}_X(D;s)$ lives, for example, in a polynomial ring `with
fractional powers' in a variable $T$ over  some ring ${\Cal R}_d$,
localized with respect to the elements $L^\nu T^N - 1$ and $L^\nu
- T^N$ for $\nu \in \frac 1 d \Bbb Z$ (and $\nu \geq 1$) and $N
\in \Bbb Q_{>0}$. (We verified that the defining expression does
not depend on the chosen resolution using the weak factorization
theorem [AKMW][W\l].)

(ii) One can specialize ${\Cal Z}_X(D;s)$ to the Hodge polynomial
level via the map $H$, obtaining
$$Z_X (D;s) = (uv)^{-n} \sum_{I \subset S} H(E^\circ_I) \prod_{i
\in I} \frac{uv-1}{(uv)^{\nu_i + sN_i}-1},$$ where now $(uv)^{-s}$
is a variable.

(iii) For any constructible subset $W$ of $X$, we can consider
more generally zeta functions ${\Cal Z}_{W \subset X} (D;s)$ and
$Z_{W \subset X}(D;s)$, using $E^\circ_I \cap h^{-1} W$ instead of
$E^\circ_I$ in the defining expressions.

\bigskip
\noindent {\it Note.} Here we re-normalized the zeta functions of
[Ve3, \S 2] with a factor $L^{-n}$ and $(uv)^{-n}$, respectively.
\bigskip
\bigskip
\noindent {\bf 2. Motivic principal value integrals on smooth
varieties}
\bigskip
\noindent {\bf 2.1.} Let $Y$ be a smooth algebraic variety of
dimension $n$. Let $\omega^{1/d}$ be a multi-valued differential
form on $Y$, such that $\div \omega^{1/d}$ is a normal crossings
divisor and $\omega^{1/d}$ has no logarithmic poles.

Denote by $E_i, i \in S$, the irreducible components of $E = \supp
(\div \omega^{1/d})$, and let $\alpha_i - 1$ be the multiplicity
of $E_i$ in $\div(\omega^{1/d})$. So $\div \omega^{1/d} = \sum_{i
\in S} (\alpha_i - 1) E_i$, and the $\alpha_i \in \frac 1 d \Bbb Z
\setminus \{0\}$. For $I \subset S$ we put $E^\circ_I := (\cap_{i
\in I} E_i) \setminus (\cup_{\ell \not\in I} E_\ell)$. Note that
$Y = \coprod_{I \subset S} E^\circ_I$.

\bigskip
\noindent {\bf 2.2.} We start by giving two equivalent definitions
for the principal value integral of $\omega^{1/d}$ on $Y$ on the
level of Hodge polynomials. The first one is analogous to the
classical real and $p$-adic situation, and the second one will
turn out to be a specialization of the definition on the motivic
level.

\bigskip
(1) Consider for $s \in \Bbb Z, s >\!> 0$, the `motivic integral
on Hodge polynomial level'
$$I(s) := \int_{{\Cal L}(Y)} (uv)^{-\ord_t (\div \omega^{1/d} +
sE)} d \mu ,$$ where ${\Cal L}(Y)$ is the arc space of $Y$ and
$\ord_t(\cdot)$ denotes the order of the given divisor along an
arc in ${\Cal L}(Y)$; see e.g. [DL2][Ba][Ve2]. This is the
analogue of the converging integral for $s >\!> 0$ in the
classical case. By [Ba] or [DL2] also $I(s)$ converges for $s >\!>
0$ (in fact if and only if $\alpha_i + s > 0$ for all $i \in S$),
and then
$$I(s) = (uv)^{-n} \sum_{I \subset S} H(E^\circ_I) \prod_{i \in I}
\frac{(uv - 1)(uv)^{-s}}{(uv)^{\alpha_i} - (uv)^{-s}}.$$ Consider
now the unique rational function $Z(T)$ over $\Bbb
Q(u^{1/d},v^{1/d})$ in the variable $T$, yielding $I(s)$ when
evaluated in $T = (uv)^{-s}$ for all $s >\!> 0$; it is given by
$$Z(T) = (uv)^{-n} \sum_{I \subset S} H(E^\circ_I) \prod_{i \in I}
\frac{(uv-1)T}{(uv)^{\alpha_i} - T}.$$

\bigskip
\noindent {\bf Hodge level definition 1.} The principal value
integral of $\omega^{1/d}$ on $Y$ is $\lim_{T \rightarrow 1} Z(T)
= ev_{T=1} Z(T)$ and is thus given by the formula
$$(uv)^{-n} \sum_{I \subset S} H(E^\circ_I) \prod_{i \in I}
\frac{uv - 1}{(uv)^{\alpha_i}-1}.$$ Note that this proces is
indeed analogous to the classical case, where we take the limit
for $s \rightarrow 0$ of a meromorphic continuation.

\bigskip
(2) We consider the zeta function $Z_Y(\div \omega^{1/d}; s)$ of
(1.6(ii)). Since here $\div \omega^{1/d}$ is already a normal
crossings divisor on the smooth variety $Y$, we have
$$Z_Y (\div \omega^{1/d};s) = (uv)^{-n} \sum_{I \subset S}
H(E^\circ_I) \prod_{i \in I} \frac{uv-1}{(uv)^{1+(\alpha_i-1)s}  -
1}.$$ (If you don't like fractional powers of $T = (uv)^{-s}$,
just consider $T^{1/d}$ as a variable with integer powers
$-d(\alpha_i - 1)$.)

\bigskip
\noindent {\bf Hodge level definition 2.} The principal value
integral of $\omega^{1/d}$ on $Y$ is $$\lim_{s \rightarrow 1} Z_Y
(\div \omega^{1/d};s) = ev_{s=1} Z_Y (\div \omega^{1/d};s).$$ This
means of course evaluating in $T = (uv)^{-1}$, and yields the same
formula as in the previous definition.

\bigskip
\noindent {\sl Remark.} Alternatively, we could have taken the
(re-normalized) zeta function of [Ve2], associated to the
effective divisor $aE$ and the sheaf of {\it regular} differential
forms ${\Cal O}(aE) \otimes \omega^{1/d}$ for some $a >\!> 0$,
which is given by the formula
$$(uv)^{-n} \sum_{I \subset S} H(E^\circ_I) \prod_{i \in I}
\frac{uv - 1}{(uv)^{a+\alpha_i+sa}-1},$$ and evaluate it in
$s=-1$.

\bigskip
\noindent {\bf 2.3.} On the level of the Grothendieck ring, we
cannot use the first approach since zero divisors may occur. The
second approach however generalizes and yields the desired
formula.

\bigskip
\noindent {\bf Definition.} The motivic principal value integral
of $\omega^{1/d}$ on $Y$ is the evaluation of ${\Cal Z}_Y(\div
\omega^{1/d};s)$
 in $s=1$; it is given by the formula
 $$L^{-n} \sum_{I \subset S} [E^\circ_I] \prod_{i \in I}
 \frac{L-1}{L^{\alpha_i}-1},$$
 living in the ring ${\Cal R}_d$ of (1.5). We denote it by $PV
 \int_Y \omega^{1/d}$.

 \bigskip
 \noindent
 {\sl Remarks.} (1) One easily verifies that evaluating ${\Cal Z}_Y
 (\div \omega^{1/d};s)$ in $s = 1$ (i.e. in $T = L^{-1}$) indeed
 yields a well defined element in ${\Cal R}_d$.

 (2) In the special case that all $\alpha_i > - 1$, we can just
 use the converging motivic integral $\int_{{\Cal L}(Y)}
 L^{-\ord_t(\div \omega^{1/d})} d \mu$, given by the same formula,
 but then it is only well defined in a completion of $K_0 (\Var_k)
 [L^{-1/d}]$, see [DL2] and [Ve2].

 (3) For a constructible subset $W$ of $Y$ we can consider more
 generally $PV \int_{W \subset Y} \omega^{1/d}$ as the evaluation
 of ${\Cal Z}_{W \subset Y}(\div \omega^{1/d};s)$ in $s=1$, see
 (1.6(iii)). This corresponds morally to `classical' $p$-adic and
 real principal value integrals involving a locally constant
 function and $C^\infty$ function, respectively, with compact
 support.

 \bigskip
 \noindent
 {\bf 2.4.} Let now $X$ be a smooth algebraic variety of dimension
 $n$, and $\omega^{1/d}$ any multi-valued differential form on $X$.
 Is it possible to associate a well defined `natural' principal
 value integral to $\omega^{1/d}$ ? Of course in some sense
 logarithmic poles will have to be excluded.

 A first natural idea is to consider (the pull-back of)
 $\omega^{1/d}$ via a modification $h : Y \rightarrow X$ such that
 the divisor of $h^\ast \omega^{1/d}$ is a normal crossings
 divisor on $Y$. If there exists a modification for which
 $h^\ast \omega^{1/d}$ has no logarithmic poles on $Y$, then
 the desired principal value integral could be defined as $PV
 \int_Y h^\ast \omega^{1/d}$. Of course the point here is whether
 this is independent of the chosen modification.

 Looking at the zeta function approach in Definition 2.3, another
 natural idea is just to define the principal value integral of
 $\omega^{1/d}$ on $X$ as the evaluation of ${\Cal Z}_X (\div
 \omega^{1/d};s)$ in $s=1$, {\it if this makes sense.} Here no
 choices are involved.

 We verify that the first approach works and yields the same
 result as the second approach {\sl if} there exists a modification $h :
 Y \rightarrow X$ satisfying a slightly stronger condition than
 above. Then we indicate some subtle problems concerning the
 `naive' first approach.

 \bigskip
 \noindent
 {\bf 2.5.} We say that a modification $h : Y \rightarrow X$ is
 {\sl good} if it is a log resolution of $D := \div \omega^{1/d}$ on
 $X$, for which $\div(h^\ast \omega^{1/d})$ has no logarithmic
 poles on $Y$.

 Note that there could exist modifications $h : Y \rightarrow X$
 for which $\div(h^\ast \omega^{1/d})$ is a normal crossings
 divisor and $h^\ast \omega^{1/d}$ has no logarithmic poles, but
 such that $h^\ast D$ is {\sl not} a normal crossings divisor ! Indeed,
 it is possible that $h^\ast \omega^{1/d}$ has multiplicity zero
 along some component of $h^\ast D$, meaning that this component
 does not occur in $\div(h^\ast \omega^{1/d})$. For example, let
 $D$ (locally) be given by $\frac 1 2 D_1 -  \frac 3 2 D_2$, where
$D_1$ and $D_2$ are smooth curves, intersecting each other in a
point $P$ with intersection multiplicity 2, see Figure 1. Let $h$
be the blowing-up of $X$ with centre $P$. Then $h^\ast
\omega^{1/d}$ has multiplicity zero along the exceptional curve
$D_3$ of $h$; so $\div(h^\ast \omega^{1/d})$ is (locally) a normal
crossings divisor but $h^\ast D$ is not.

\vskip 1truecm \centerline{\beginpicture
 \setcoordinatesystem units <.5truecm,.5truecm> point at -10 0
 \setquadratic   \plot 1 1.5   0 0   1 -1.5 /
                \plot -1 1.5   0 0   -1 -1.5 /
 \put {$\bullet$} at 0 0
 \put {$P$} at  .6 0
 \put{$D_1$} at -1.6 1.5
 \put{$D_2$} at 1.7 1.5
 \put{$\longleftarrow$} at 5 0
 \put{$h$} at 5 .6
 \setcoordinatesystem units <.5truecm,.5truecm> point at -21 0
 \setdashes
 \putrule from -1.9 0 to 1.9 0
 \setsolid
 \setlinear \plot -1 -1.5   1 1.5 /
 \plot 1 -1.5   -1 1.5 /
  \put{$D_1$} at -1.6 1.5
 \put{$D_2$} at 1.7 1.5
 \put{$D_3$} at 2.6 0
  \put {Figure 1} at -6 -2.5
 \endpicture}

\vskip 1truecm

\proclaim{2.6. Proposition} Suppose that there exists a good
modification $h : Y  \rightarrow X$. Denote by $E_i, i \in S$, the
irreducible components of $h^{-1}D$, and let $\alpha_i - 1$ be the
multiplicity of $E_i$ in $\div(h^\ast \omega^{1/d})$. Put
$E^\circ_I := (\cap_{i \in I} E_i) \setminus (\cup_{\ell \not\in
I} E_\ell)$ for $I \subset S$. Then the evaluation of ${\Cal Z}_X
(\div \omega^{1/d};s)$ in $s=1$ is well defined, it is equal to
$PV \int_Y h^\ast \omega^{1/d}$, and is given by the formula
$L^{-n} \sum_{I \subset S} [E^\circ_I] \prod_{i \in I}
\frac{L-1}{L^{\alpha_i}-1}$.
\endproclaim
\medskip \demo{Note} Here $\alpha_i = 1$ could occur,
meaning that $E_i$ does not appear in $\div(h^\ast \omega^{1/d})$.
\enddemo

\medskip
 \demo{Proof} Denote the multiplicities of $E_i$ in
$K_{Y|X}$ and in $h^\ast D$ by $\nu_i - 1$ and $N_i$,
respectively. Since $h$ is really  a log resolution of $D$, we can
express ${\Cal Z}_X(D;s)$ as
$$L^{-n} \sum_{I \subset S} [E^\circ_I] \prod_{i \in I}
\frac{L-1}{L^{\nu_i + sN_i}-1}.$$ Now $D = \div \omega^{1/d}$ and
$\div(h^\ast \omega^{1/d})$ are representatives of $K_X$ and
$K_Y$, respectively. Hence $\div(h^\ast \omega^{1/d}) = K_{Y|X} +
h^\ast D$, meaning that $\alpha_i = \nu_i  + N_i$ for all $i \in
S$. Since  $\div(h^\ast \omega^{1/d})$ has no logarithmic poles,
all $\alpha_i \ne 0$. So indeed evaluating $\Cal Z_X (D;s)$ in
$s=1$ makes sense and yields the stated formula, which is just $PV
\int_Y h^\ast \omega^{1/d}$.  \qed
\enddemo

\bigskip
We define the principal value integral of $\omega^{1/d}$ on $X$ as
given by Proposition 2.6. For completeness we recall all data.

\bigskip
\noindent {\bf 2.7. Definition.} Let $X$ be a smooth algebraic
variety of dimension $n$ and $\omega^{1/d}$ a multi-valued
differential form on $X$. Suppose that there exists a log
resolution $h : Y \rightarrow X$ of $\div \omega^{1/d}$ on $X$,
for which $h^\ast \omega^{1/d}$ has no logarithmic poles on $Y$.
Then the principal value integral of $\omega^{1/d}$ on $X$,
denoted $PV \int_{X} \omega^{1/d}$, is given by one of the
equivalent expressions in Proposition 2.6.

\bigskip
\noindent {\sl Remarks.} (1) Also here we could proceed
alternatively using the zeta functions of [Ve2].

(2) We can proceed more generally, involving a constructible
subset $W$ of $X$, just as in Remark (3) after Definition 2.3.

(3) For real principal value integrals, Jacobs gave a similar
definition [Ja1, \S 7].

\bigskip
\noindent {\bf 2.8.} We return to the first approach. Suppose now
that there exists a modification   \linebreak
 $g : Z \rightarrow X$
such that $\div(g^\ast \omega^{1/d})$ is a normal crossings
divisor and $g^\ast \omega^{1/d}$ has no logarithmic poles on $Z$.
(Recall that then $g$ is not necessarily good.) Two subtle
questions impose themselves here.

\bigskip
I) Suppose that there exists at least one good modification of
$X$; so $PV \int_X \omega^{1/d}$ is defined. Is this principle
value integral then equal to the obvious formula associated to
$g^\ast \omega^{1/d}$ on $Z$ ? I.e., with $\div(g^\ast
\omega^{1/d}) = \sum_{i \in S_Z} (\alpha_i - 1)E_i$ and
$E^\circ_I$ for $I \subset S_Z$ as usual, is
$$PV \int_X \omega^{1/d} = L^{-n} \sum_{I \subset S_Z} [E^\circ_I]
\prod_{i \in I} \frac{L-1}{L^{\alpha_i}-1}\, ?$$
 If there exists a
modification $\pi : Y \rightarrow Z$ such that $h := g \circ \pi :
Y \rightarrow X$ is good, then the answer is yes. This can be
verified using the zeta function ${\Cal Z}_Z(\div (g^\ast
\omega^{1/d});s)$, for which we use the defining expressions on
both $Z$ and $Y$. One easily computes that evaluating in $s=1$
yields the right and left hand sides above, respectively. The
point is that `deleting irreducible components of $h^\ast(\div
\omega^{1/d})$ on $Y$ with $\alpha = 1$' does not change the
formula for $PV \int_X \omega^{1/d}$ in terms of $h$.

When there does not exist such a modification $Y \rightarrow Z$,
we do not know the answer.

\bigskip
II) Suppose on the other hand that {\sl no} good modification of
$X$ exists. Are the expressions
$$L^{-n} \sum_{I \subset S_Z} [E^\circ_I] \prod_{i \in I}
\frac{L-1}{L^{\alpha_i}-1}$$ for modifications $Z \rightarrow X$
as above {\it independent} of the chosen $Z$ ?

\bigskip
\noindent {\bf 2.9.} The principal value integrals of Definition
2.3 and the more general Definition 2.7 satisfy a `Poincar\'e
duality'. Bittner [Bi] showed that there exists a ring involution
${\Cal D}$ of $K_0 (\Var_k)[L^{-1}]$ satisfying ${\Cal D}(L) =
L^{-1}$, and characterized by ${\Cal D}([X]) = L^{-\dim X}[X]$
when $X$ is a complete connected smooth variety. (It is in fact a
lifting of the usual duality operator on the level of motives to
the level of varieties. And specializing ${\Cal D}$ to the level
of Hodge polynomials is just a reformulation of Poincar\'e and
Serre duality.) It extends to the rings ${\Cal R}_d$ via ${\Cal
D}(L^{1/d}) = L^{-1/d}$.

\proclaim{Proposition} The principal value integrals of
Definitions 2.3 and 2.7 satisfy $${\Cal D}(PV \int_X \omega^{1/d})
= L^{-\dim X} PV \int_X \omega^{1/d}.$$
\endproclaim

\medskip \demo{Proof} This follows from the concrete formula
for $PV \int_X \omega^{1/d}$ by the same computation as in e.g.
[Ba], [DM] or [Ve3]. \qed
\enddemo

\pagebreak

\noindent {\bf 3. Birational invariance ?}
\bigskip

Here we assume $k$ to be algebraically closed.

\bigskip
 \noindent {\bf 3.1.}
Actually, a (multi-valued) differential form is a birational
notion. When we consider such a form $\omega^{1/d}$ on a variety
$X$ and its pull-back $h^\ast \omega^{1/d}$ on $Y$ via a
modification   \linebreak
 $h : Y \rightarrow X$, this is just a
matter of notation : $\omega^{1/d}$ and $h^\ast \omega^{1/d}$ are
in fact the same element in a one-dimensional vector space over
the function field of $X$.

Fix a form $\omega^{1/d}$ for which there exists a smooth complete
variety $X$ such that $\div (\omega^{1/d})$ on $X$ is a normal
crossings divisor and $\omega^{1/d}$ has no logarithmic poles on
$X$. Then we can consider $PV \int_X \omega^{1/d}$, and it is a
natural question whether this notion depends on the chosen such
model $X$. In other words : is the motivic principal value
integral a birational invariant ?

\bigskip
\noindent {\bf 3.2.} {\sl Remark.} A necessary condition is of
course that, if $\pi : X^\prime \rightarrow X$ is the blowing-up
of an $X$ as  above in a smooth centre that has normal crossings
with $\div (\omega^{1/d})$ and such that $\omega^{1/d}$ has also
no logarithmic poles on $X^\prime$, then $PV \int_X \omega^{1/d} =
PV \int_{X^\prime} \omega^{1/d}$. This can be verified by
straightforward computations as in [Ve1], [Ve3, Lemma 2.3.2] or
[Al].

We should remark that this is however not sufficient to derive
birational invariance with the help of the weak factorization
theorem. Indeed, on some `intermediate' varieties connecting two
such models the form $\omega^{1/d}$ could have logarithmic poles.

\bigskip
\noindent {\bf 3.3.} Note that in dimension one there is only one
smooth complete model in  a given birational equivalence class. So
from now on we work in dimension at least two.

\bigskip
First we show that when the Kodaira dimension is $-\infty$, the
answer is in general negative.

\bigskip
\noindent {\bf 3.4.} {\sl Example.} We work in the class of
rational surfaces and take $\omega^{1/2}$ on $\Bbb P^2_{(X:Y:Z)}$,
given by $\omega^{1/2} = y ^{-3/2} dx dy$ on the affine chart
$\Bbb A^2_{(x,y)}$. Then on the chart $\Bbb A^2_{(y,z)}$ we have
that $\omega^{1/2} = y^{-3/2} z^{-3/2} dydz$. Denoting $C_1 := \{
Y = 0 \}$  and $C_2 := \{ Z = 0 \}$, we see that
$\div(\omega^{1/2}) = - \frac 32 C_1 - \frac 32 C_2$ is a normal
crossings divisor on $\Bbb P^2$ and  that no logarithmic poles
occur.

Consider now the birational map $\pi : \Bbb
P^2_{(X^\prime:Y^\prime:Z^\prime)} -\! \rightarrow \Bbb
P^2_{(X:Y:Z)}$ given by $\Bbb A^2_{(x^\prime,y^\prime)}
\rightarrow \Bbb A^2_{(x,y)} : (x^\prime, y^\prime) \mapsto
(x^\prime,y^\prime - x^{\prime 2})$. On $\Bbb
P^2_{(X^\prime:Y^\prime:Z^\prime)}$ our form $\omega^{1/2}$ is
given by $(y^\prime - x^{\prime 2})^{-3/2} dx^\prime dy^\prime$
and $(y^\prime z^\prime - 1)^{-3/2} dy^\prime dz^\prime$ on the
analogous charts. Hence on the `new' $\Bbb P^2$ we have that
$\div(\omega^{1/2}) = - \frac 32 C_1$, where we denote the
birational transform of $C_1$, i.e. $\{ Y^\prime Z^\prime -
X^{\prime 2} = 0 \}$, by the same symbol. The formula of
Definition 2.3 yields

$$PV \int_{\Bbb P^2_{(X:Y:Z)}} \omega^{1/2} = L^{-2} \big(L^2 - L + 2L
\frac{L-1}{L^{-1/2}-1} + \frac{(L-1)^2}{(L^{-1/2}-1)^2}\big) = 0$$
and

$$PV \int_{\Bbb P^2_{(X^\prime:Y^\prime:Z^\prime)}} \omega^{1/2} =
L^{-2} \big(L^2 + (L + 1) \frac{L-1}{L^{-1/2}-1}\big) = - L^{-3/2}
(L + L^{1/2} + 1) \ne 0.
$$

\vskip 1truecm
 \centerline{\beginpicture
 \setcoordinatesystem units <.43truecm,.43truecm>
 \putrectangle corners at 0 0 and 9 6
 \ellipticalarc axes ratio 3:2 360 degrees from 4.5 5 center at 4.5 3
 \putrule from 1 5 to 8 5
 \put {$\bullet$} at 4.5 5
  \put {$C_1$} at 2 1.2
  \put {$C_5$} at 1 4.4
  \put {$P^2_{(X':Y':Z')}$} at 11.4 5.5
  \arrow <.3truecm> [.2,.6] from 4.5 10 to 4.5 7
  \put {$ - \rightarrow$} at 16.5 2.5
  \put {$\pi$} at 16.5 3.1
  \put {Figure 2} at 16.5 -2
 \setcoordinatesystem units <.43truecm,.43truecm> point at -24 0
 \putrectangle corners at 0 0 and 9 6
 \putrule from 1 2 to 8 2
 \putrule from 7 1 to 7 5
 \put {$\bullet$} at 7 4
  \put {$C_1$} at 1.5 1.4
  \put {$C_2$} at 7.7 3
  \put {$P^2_{(X:Y:Z)}$} at -1.8 5.5
  \arrow <.3truecm> [.2,.6] from 4.5 10 to 4.5 7
 \setcoordinatesystem units <.43truecm,.43truecm> point at 0 -11
 \putrectangle corners at 0 0 and 9 7
 \putrule from 1 4 to 8 4
 \putrule from 6 1 to 6 6
 \setlinear  \plot 1 1  8 5.2 /
 \put {$\bullet$} at 6 4
  \put {$C_1$} at 1.5 2.1
  \put {$C_3$} at 1.5 4.6
  \put {$C_5$} at 6.7 1
  \put {$S_5$} at 9.7 6.5
  \arrow <.3truecm> [.2,.6] from 8 11 to 6 8
 \setcoordinatesystem units <.43truecm,.43truecm> point at -12 -11
 \putrectangle corners at 0 0 and 9 7
 \putrule from 1 2.6 to 8 2.6
 \putrule from 1 5.4 to 8 5.4
 \putrule from 6 1 to 6 6
  \multiput {$\bullet$} at 6 2.6  6 1.8 /
  \put {$C_1$} at 1.5 2
  \put {$C_3$} at 1.5 6
  \put {$C_4$} at 6.7 4
  \put {$S_3$} at 9.7 6.5
  \arrow <.3truecm> [.2,.6] from 8 11 to 6 8
  \arrow <.3truecm> [.2,.6] from 1 11 to 3 8
 \setcoordinatesystem units <.43truecm,.43truecm> point at -24 -11
 \putrectangle corners at 0 0 and 9 7
 \putrule from 1 2 to 8 2
 \putrule from 1 5 to 8 5
 \putrule from 6 1 to 6 6
  \put {$\bullet$} at 6 5
  \put {$C_1$} at 1.5 1.4
  \put {$C_3$} at 1.5 5.5
  \put {$C_2$} at 6.8 3.5
  \put {$S_1$} at 9.7 6.5
  \arrow <.3truecm> [.2,.6] from 1 11 to 3 8
 \setcoordinatesystem units <.43truecm,.43truecm> point at -6 -23
 \putrectangle corners at 0 0 and 9 8
 \putrule from 1 5.8 to 8 5.8
 \putrule from 6 1 to 6 7
 \setlinear  \plot 1 3  8 4.2 /
      \plot 3 1   8 2.4 /
  \put {$C_1$} at 1.5 3.7
  \put {$C_3$} at 1.5 6.3
  \put {$C_4$} at 6.8 6.8
  \put {$C_5$} at 2.4 1
  \put {$S_4$} at 9.7 7.5
 \setcoordinatesystem units <.43truecm,.43truecm> point at -18 -23
 \putrectangle corners at 0 0 and 9 8
  \setlinear \plot 1 2  7 1 /  \plot  1 6  7 7 /
     \plot 5.5 .6  8 5 /
      \plot 5.5 7.4   8 3 /
  \put {$C_1$} at 1.5 1.3
  \put {$C_3$} at 1.5 6.7
  \put {$C_4$} at 6.1 5.3
  \put {$C_2$} at 6 2.9
  \put {$S_2$} at 9.7 7.5
 \endpicture}
\vskip .8truecm

It is useful to indicate the `geometric reason' why this happens.
We decompose in Figure 2 the map $\pi$ in  a composition of
blowing-ups and blowing-downs, where the fat points indicate the
centers of blowing-up and $C_3,C_4$ and $C_5$ are  exceptional
curves. The surfaces on the middle row are ruled surfaces. One
easily verifies that the multiplicities of $C_3,C_4$ and $C_5$ in
$\div \omega^{1/2}$ are $- \frac 12, - 1$ and $0$, respectively.
This means in particular that $C_5$ does not occur in the support
of $\div \omega^{1/2}$ on $S_4, S_5$ and $\Bbb
P^2_{(X^\prime:Y^\prime:Z^\prime)}$ (where $C_5$ is a line), and
that $\div \omega^{1/2}$ has logarithmic poles on $S_2,S_3$ and
$S_4$. So we {\sl cannot} identify $PV \int_{\Bbb P^2_{(X:Y:Z)}}
\omega^{1/2} = PV \int_{S_1} \omega^{1/2}$ with $PV \int_{\Bbb
P^2_{(X^\prime:Y^\prime:Z^\prime)}} \omega^{1/2} = PV \int_{S_5}
\omega^{1/2}$ via principal value integrals on the intermediate
surfaces $S_2,S_3$ and $S_4$ since they are {\sl not defined}
there.

\medskip
 \noindent {\sl Note.} (i) We found this example several
years ago; it was briefly mentioned by Jacobs [Ja1, \S 8] in the
context of real principal value integrals.

(ii) Actually, it is also valid when $k$ is not algebraically
closed. (And if we would have introduced principal value integrals
in arbitrary characteristic by the same formula, it would still
work.)


\bigskip
\noindent {\bf 3.5.} Example 3.4 can be adapted to the birational
equivalence class of any non-rational ruled surface; then we use
only the middle and top row of Figure 2. It is possible to give a
similar form $\omega^{1/2}$ with $PV \int_{S_1} \omega^{1/2} = 0$
and $PV \int_{S_5} \omega^{1/2} \ne 0$ (maybe $\div(\omega^{1/2})$
will contain more fibres in its support). So for surfaces such
examples exist in every birational equivalence class of Kodaira
dimension $-\infty$.

Moreover, by taking Cartesian products with arbitrary complete
smooth varieties, Example 3.4 can be extended to arbitrary
dimension.

\bigskip
We now turn to the other case, i.e. when the Kodaira dimension is
nonnegative.

\bigskip
\proclaim {3.6. Theorem} Fix the birational equivalence class of a
{\it surface} of nonnegative Kodaira dimension, and a multivalued
differential form $\omega^{1/d}$ on it. Suppose that this class
contains a smooth complete model $X$ such that
$\div(\omega^{1/d})$ is a normal crossings divisor on $X$ and
$\omega^{1/d}$ has no logarithmic poles on $X$. Then $PV \int_{X}
\omega^{1/d}$ is independent of the chosen such model, and is thus
a birational invariant.
\endproclaim

\medskip \demo{Proof} (Recall that for smooth surfaces complete
is equivalent to projective.) Let $X_m$ be the unique (smooth,
projective) minimal model in the class, and let $h : Y \rightarrow
X_m$ be the composition of the minimal set of blowing-ups, needed
to make $\div \omega^{1/d}$ a normal crossings divisor on $Y$.
More precisely, if $\div \omega^{1/d}$ is a normal crossings
divisor on $X_m$, put $Y := X_m$. Otherwise, let $h_1 : Y_1
\rightarrow X_m$ be obtained from $X_m$ by blowing up the finite
number of points where $\div \omega^{1/d}$ has no normal
crossings. If $\div \omega^{1/d}$ is a normal crossings divisor on
$Y_1$, put $Y := Y_1$. Otherwise, continuing this way abuts in the
unique smooth projective $Y$ such that $\div \omega^{1/d}$ is a
normal crossings divisor on $Y$ and $h : Y \rightarrow X_m$ is
minimal with respect to this property.

Suppose now that $X$ is any model as in the \'enonc\'e of the
theorem. Since $\div \omega^{1/d}$ is a normal crossings divisor
on $X$, there is a morphism $\pi : X \rightarrow Y$. Also, since
$\omega^{1/d}$ has no logarithmic poles on $X$, the same is
certainly  true on $Y$. Hence also $PV \int_Y \omega^{1/d}$ is
defined, and is equal to $PV \int_X \omega^{1/d}$ by Remark 3.2.
\qed
\enddemo

\bigskip
\noindent {\bf 3.7.} In higher dimensions we face the
non-existence of a unique minimal model, and the fact that in
general a (smooth, complete) variety does not map to a minimal
model by a {\it morphism}. A reasonable idea is to try to adapt
Definition 2.7 as follows, assuming the Minimal Model Program.
Take a birational equivalence class of nonnegative Kodaira
dimension and a multi-valued differential form $\omega^{1/d}$ on
it.

Suppose that there exists a minimal model $X$ in the given class,
and a log resolution $h : Y \rightarrow X$ of $\div \omega^{1/d}$
on $X$, for which $\omega^{1/d}$ has no logarithmic poles on $Y$.
Define then the (candidate) birational invariant associated to
$\omega^{1/d}$ as $PV \int_Y \omega^{1/d}$. Now independence of
both $X$ and $Y$ has to be checked.

Fixing such a minimal model $X$, the independence of $Y$ is proven
by the same argument as for Proposition 2.6. Indeed, $PV \int_Y
\omega^{1/d}$ is just the evaluation of ${\Cal Z}_X (\div
\omega^{1/d};s)$ in $s=1$. (Recall that this zeta function was
defined more generally on canonical varieties, so certainly on
minimal models.) We now verify that the zeta function ${\Cal Z}_X
(\div \omega^{1/d};s)$ itself in fact does not depend on the
chosen $X$; then a fortiori the same is true for its evaluation in
$s=1$.

\bigskip
\noindent \proclaim {\bf 3.8. Proposition} Let $\omega^{1/d}$ be a
multi-valued differential form on  a birational equivalence class
of nonnegative Kodaira dimension. Let $X_1$ and $X_2$ be minimal
models in this class; then ${\Cal Z}_{X_1} (\div \omega^{1/d};s) =
{\Cal Z}_{X_2} (\div \omega^{1/d};s)$.
\endproclaim

\medskip
\demo{Proof} Since $X_1$ and $X_2$ are isomorphic in codimension
one, $\div \omega^{1/d}$ on $X_1$ and $X_2$ are each others
birational transform. Take a common log resolution $h_i : Y
\rightarrow X_i$ of $\div \omega^{1/d}$ on $X_1$ and $X_2$.
Remembering that these two $(\Bbb Q$-)divisors are representatives
of $K_{X_1}$ and $K_{X_2}$, we have $h^\ast_1 (\div \omega^{1/d})
= h^\ast_2 (\div \omega^{1/d})$ and $K_{Y|X_1} = K_{Y|X_2}$, see
[KM, Proof of Theorem 3.52] or [Wa, Corollary 1.10].

So, computing both zeta functions by their defining expressions on
$Y$, see (1.6), yields exactly the same formula. \qed
\enddemo

\bigskip
\noindent {\sl Note.} When our birational equivalence class is of
general type, we can use its {\it unique} canonical model $X_c$
and work with ${\Cal Z}_{X_c} (\div \omega^{1/d};s)$, which
shortens the argument.

\bigskip
Summarizing, we obtained the following well defined invariant.

\bigskip
\noindent {\bf 3.9. Definition.} Let $\omega^{1/d}$ be a
multi-valued differential form on a birational equivalence class
of nonnegative Kodaira dimension. Assume the Minimal Model
Program. Suppose that there exists a minimal model $X$ and a log
resolution $h : Y \rightarrow X$ of $\div \omega^{1/d}$ on $X$,
for which $\omega^{1/d}$ has no logarithmic poles on $Y$. Then $PV
\int_Y \omega^{1/d}$ is independent of all choices and is thus a
(partial) birational invariant of $\omega^{1/d}$.

\bigskip
\noindent {\bf 3.10.} We mention `partial' because we are
confronted with similar subtle problems as in $\S 2$. Suppose that
there exists a smooth complete model $Z$ in the given class, for
which $\div \omega^{1/d}$ is a normal crossings divisor on $Z$ and
$\omega^{1/d}$ has no logarithmic poles on $Z$.

I) If there exists a minimal model $X$ and a log resolution $Y
\rightarrow X$ as in Definition 3.9, is the invariant above then
equal to $PV \int_{Z} \omega^{1/d}$ ?

II) If on the other hand no such $X$ and $Y$ exist, are the
expressions $PV \int_Z \omega^{1/d}$ the same for different such
models $Z$ ?

\bigskip
\noindent {\bf 3.11.} {\sl Remark.} An alternative point of view
for `birational invariance' is `independence of chosen completion'
for principal value integrals on non-complete smooth varieties.
For real principal value integrals Jacobs [Ja1, \S 8] mentioned
Example 3.4 in this context.


\bigskip
\bigskip
\Refs

\ref \key AKMW \by D. Abramovich, K. Karu, K. Matsuki, J.
Wlodarczyk \paper Torification and factorization of birational
maps \jour J. Amer. Math. Soc. \vol 15 \yr 2002 \pages 531--572
\endref

\ref\key Al
 \by P\. Aluffi
 \paper Chern classes of birational varieties
 \jour math.AG/0401167
 \endref

\ref\key AVG \by V\. Arnold, A\. Varchenko and S\.
Goussein--Zad\'e \book Singularit\'es des applications
diff\'erentiables II \publ Editions Mir \publaddr Moscou \yr 1986
\endref

\ref \key Ba \by V\. Batyrev \paper Stringy Hodge numbers of
varieties with Gorenstein canonical singularities \jour Proc.
Taniguchi Symposium 1997, In \lq Integrable Systems and Algebraic
Geometry, Kobe/Kyoto 1997\rq, World Sci. Publ. \vol \yr 1999
\pages 1--32
\endref

 \ref \key Bi
 \by F. Bittner
 \paper The universal Euler
characteristic for varieties of  characteristic zero
 \jour Compositio Math.
 \vol
 \yr to appear
 \pages
\endref

\ref \key Cr
 \by A\. Craw
 \paper An introduction to motivic
integration \jour math.AG/9911179 \vol \yr 2001 \pages
\endref

\ref \key De1 \by J\. Denef \paper On the degree of Igusa's local
zeta function \jour Amer. J. Math. \vol 109 \yr 1987 \pages
991--1008
\endref

\ref \key De2 \by J\. Denef \paper Report on Igusa's local zeta
function \jour Ast\'erisque \paperinfo S\'em. Bourbaki 741 \vol
201/202/203 \yr 1991 \pages 359--386
\endref

\ref \key DJ
 \by J\. Denef and Ph\. Jacobs
 \paper On the vanishing of principal value integrals
 \jour C. R. Acad. Sci. Paris
 \vol 326
 \yr 1998
 \pages 1041--1046
 \endref

\ref \key DL1 \by J\. Denef and F\. Loeser \paper Motivic Igusa
zeta functions \jour J. Alg. Geom. \vol 7 \yr 1998 \pages 505--537
\endref

\ref \key DL2 \by J\. Denef and F\. Loeser \paper Germs of arcs on
singular algebraic varieties and motivic integration \jour Invent.
Math. \vol 135 \yr 1999 \pages 201--232
\endref


\ref \key DL3 \by J\. Denef and F\. Loeser \paper Geometry on arc
spaces of algebraic varieties \paperinfo Proceedings of the Third
European Congress of Mathematics, Barcelona 2000 \jour Progr.
Math. \vol 201 \publ Birkh\"auser, Basel \yr 2001 \pages 327--348
\endref

\ref \key DM \by J\. Denef and D\. Meuser \paper A functional
equation of Igusa's local zeta function \jour Amer. J. Math. \vol
113 \yr 1991 \pages 1135--1152
\endref


\ref \key Ha1
 \by T\. Hales
 \paper Can p-adic integrals be computed?
 \jour "Contributions to Automorphic Forms, Geometry and Arithmetic"
dedicated to J. Shalika
 \publ Johns Hopkins University Press (math.RT/0205207)
 \yr
 \endref

\ref \key Ha2
 \by T\. Hales
 \paper Orbital Integrals are Motivic
 \jour math.RT/0212236
 \endref


\ref \key Hi \by H\. Hironaka \paper Resolution of singularities
of an algebraic variety over a field of
       characteristic zero
\jour Ann. Math. \vol 79 \yr 1964 \pages 109--326
\endref

\ref \key Ig1 \by J\. Igusa \paper Complex powers and asymptotic
expansions I \jour J. Reine Angew. Math. \vol 268/269 \yr 1974
\pages 110--130 \moreref \paper II \jour ibid. \vol 278/279 \yr
1975 \pages 307--321
\endref

\ref\key Ig2 \by J\. Igusa \paper Lectures on forms of higher
degree \jour Tata Inst. Fund. Research, Bombay \vol \yr 1978
\pages
\endref

\ref \key Ja1
 \by Ph\. Jacobs
 \paper Real principal value integrals
 \jour Monatsch. Math.
 \vol 130
 \yr 2000
 \pages 261--280
 \endref

 \ref \key Ja2
 \by Ph\. Jacobs
 \paper The distribution $|f|^\lambda$, oscillating integrals and
 principal value integrals
 \jour J. Analyse Math.
 \vol 81
 \yr 2000
 \pages 343--372
 \endref


\ref \key KM \by J\. Koll\'ar and S\. Mori \book Birational
geometry of algebraic varieties \bookinfo Cambridge Tracts in
Mathematics 134 \publ Cambridge Univ. Press \yr 1998
\endref


\ref \key Ko \by M\. Kontsevich \paper \jour Lecture at Orsay
(December 7, 1995) \yr \pages
\endref

\ref\key Lae
 \by A\. Laeremans
 \book The distribution $|f|^s$, topological zeta functions and Newton
      polyhedra
 \bookinfo Ph. D. thesis, Univ. Leuven
 \yr 1997
 \endref


\ref \key Lanl
 \by R\. Langlands
 \paper Orbital integrals on forms of $SL(3)$, I
 \jour Amer. J. Math.
 \vol 105
 \yr 1983
 \pages 465--506
 \endref

 \ref \key Lan2
 \by R\. Langlands
 \paper Remarks on Igusa theory and real orbital integrals
 \inbook  The Zeta Functions of Picard Modular Surfaces
 \publ Les Publications CRM, Montr\'eal; distributed by AMS
  \yr 1992
 \pages 335--347
 \endref

\ref \key Lo \by E. Looijenga \paper Motivic measures \jour
S\'eminaire Bourbaki \vol 874 \yr 2000 \pages
\endref

\ref \key LS1
 \by R\. Langlands and D\. Shelstad
 \paper On principal values on $p$--adic manifolds
 \jour Lect. Notes  Math.
 \vol 1041
 \publ Springer, Berlin
 \yr 1984
  \endref

\ref \key LS2
 \by R\. Langlands and D\. Shelstad
 \paper Orbital integrals on forms of $SL(3)$, II
 \jour Can. J. Math.
 \vol 41
 \yr 1989
 \pages 480--507
 \endref

\ref \key Ma \by K\. Matsuki \book Introduction to the Mori
Program \bookinfo Universitext \publ Springer-Verlag, New York \yr
2002
\endref


\ref \key Po
 \by  B\. Poonen
 \paper The Grothendieck ring of
varieties is not a domain
 \jour Math. Res. Letters
 \vol 9
 \yr 2002
 \pages 493--498
\endref

\ref \key Ve1 \by W\. Veys \paper Poles of Igusa's local zeta
function and monodromy \jour Bull. Soc. Math. France \vol 121 \yr
1993 \pages 545--598
\endref

\ref \key Ve2 \by W\. Veys \paper Zeta functions and \lq
Kontsevich invariants\rq\ on singular varieties \jour Canadian J.
Math. \vol 53 \yr 2001 \pages 834--865
\endref

\ref \key Ve3
 \by W\. Veys
 \paper Stringy zeta functions of $\Bbb Q$--Gorenstein varieties
\jour Duke Math. J. \vol 120 \yr 2003 \pages 469--514
\endref

\ref \key Ve4
 \by W\. Veys
 \paper Arc spaces, motivic integration and stringy invariants
\jour math.AG/0401374
 \vol
 \yr
  \pages
\endref

\ref \key Wa \by C.-L. Wang \paper On the topology of birational
minimal models \jour J. Differential Geom. \vol 50 \yr 1998 \pages
129--146
\endref

\ref \key W{\l}
 \by J. W{\l}odarczyk
 \paper Combinatorial
structures on toroidal varieties and a proof of the weak
factorization theorem \jour Invent. Math. \vol 154 \yr 2003 \pages
223--331
\endref

\endRefs
\enddocument